\documentclass[11pt]{article}
\usepackage{amsmath}
\usepackage{amsfonts}
\usepackage{amssymb}
\usepackage{amsthm}

\def\F{{\cal F}}
\def\P{\cal P}

\def\ZZ{\mathbb Z}
\def\QQ{\mathbb Q}
\def\CC{\mathbf C}

\newtheorem{theorem}{Theorem}

\newtheorem{lemma}{Lemma}
\newtheorem{remark}{Remark}
\newtheorem{proposition}{Proposition}

\title{Directed lattice paths avoiding periodic subset of points on ``time''-axis}
\author{S. Tarasov\thanks{{\tt e-mail: serge99meister@gmail.com}}}

\begin{document}
\maketitle

\begin{abstract} We compute generating functions of the set of directed lattice paths 
starting from the origin and

 avoiding a periodic set of even point on  $OX=\mbox{
 ``time''-axis}$.
As an application 
we prove a combinatorial identity proposed by P. Hajnal  and G.V. Nagy. 


\end{abstract}

\section{Introduction}

A {\em lattice walk} on $\ZZ^d,\,d\geq 1$ is a path starting from some lattice point and passing with steps that belong to some finite set $S \subset \ZZ^d$.
 Respectively, a {\em (directed) lattice path (DLP)} on an integer $(d+1)$-dimensional lattice $\ZZ\times \ZZ^{d+1}$ 
is a path 
starting from some lattice point
and passing with steps that belong to some finite set
 $\Sigma=\{(1,q)\},\  q\in S\subset \ZZ^d$. 
 The number of steps in a lattice walk or in DLP is called its {\em length}.
 
 It is convenient to assume that DLP proceeds on
 a direct product $\ZZ_+(T-\mbox{`'time'' axis})\times \ZZ^d (\mbox{space})$ so that
 DLP may be regarded as a ''world line'' or a space-time diagram of a lattice walk in $\ZZ^d$ with steps in a set $S$. And conversely,  projectlng DLP along $T$-axis onto ``space''-plane $\ZZ^d$ gives rise to the corresponding {\em projected} lattice walk in $\ZZ^d$ that starts from some lattice point  and uses  steps in the set $S$. 
 
By definition two lattice walks or  two DLP's are equal iff the corresponding lists consisting of a starting point and the sequence of consecutive lattice steps from $S$ or, respectively, $\Sigma$ are equal.

{\em Let us assume hereinafter that $S=\{-1,1\}^d$.}
 
In this paper we enumerate  DLP's 
in $\ZZ_+(=T-\mbox{axis})\times \ZZ^d (\mbox{space})$ 
that avoid particular points 
on $T$-axis. Let us call such a lattice and DLP's on it {\em restricted} and, respectively, let us call {\em unrestricted} 
the corresponding lattice $\ZZ^{d+1}$ and 
the corresponding DLP's with the same step set and initial point that may additionally use the forbidden points. 

Though making some lattice points forbidden  is more or less a standard subject, an analogous limitation on the ``time''-axis seems not to be a case. 

It will be convenient to specify not the set of forbidden points but their complement --- the set $A\subset \ZZ_+$ of admissible points.  We represent $A$ by an increasing sequence $\{a_0<a_1<\dots\}$, where 
$\{2a_0, 2a_1,\dots\}$ is a
 corresponding set of admissible points on $T$-axis. 


{\em Let us  assume hereinafter that $\pmod{p}$ operation returns value in the interval $(0,1,\dots,p-1)$}.

A set $A\subset \ZZ_+$ is called {\em periodic} if it can be represented as a union of even residue classes\footnote{Note that in periodic case we may assume w.l.o.g. that $a_0=0$.} $A=\{2a_0\pmod{2t_A}=0,\dots,2a_k=\pmod{2t_A}\}\subset \ZZ_+$ for some 
integer $t_A> a_k\geq 0$. With slight abuse of notation we denote a periodic set by a pair $\mathbf{A}\stackrel{def}{=}(A=\{a_0,\dots,a_k\}, t_A)$ and mark it in bold.  

Denote by $\P(\mathbf{A})$ the set of all starting from the origin DLP's of even length that may reach $T$-axis at points of $A$ only.
For instance, let $d=1$ then for the periodic set $\mathbf{A_1}=(\{0\}, 2)$ the set $\P(\mathbf{A_1})$ consists of of all 
DLP's of even length in the plane starting from the origin  that may reach $T$-axis at points whose abscissa is divisible by four only.  Respectively, for $\mathbf{A_2}=(\{0,1\}, 4)$ the set $\P(\mathbf{A_2})$ consists of of all 
DLP's of even length in the plane starting from the origin  that may reach $T$-axis at points whose abscissa belongs to the disjoint union of arithmetic progressions $\{8k,\,k=0,1,\dots\}\sqcup\{2+8l,\, l=0,1,\dots\}$.

Informally, periodicity limitation means that the corresponding projected along $T$-axis lattice walk on $\ZZ^d$ starting from the origin and using prescribed set of steps in $S$ is performed with universal  clock frequency such that at some (even) periodic moments of time it is prohibited to return to the origin (form a current loop passing through the origin of $\ZZ^d$). 

The paper is organized as follows. After introduction in section \ref{1} we reduce DLP's enumeration to solving an infinite in general system of linear equations in generating functions in the ring $\QQ[\![t]\!]$. In particular, if the set of forbidden points on $T$-axis is periodic then the system involved is finite and the system's matrix is a circulant matrix and have a  unique solution. Moreover we show how to obtain a set of the consequences of the system of the same dimension in which all the coefficients of the system's variables are particular multisections of the generating functions involved. In the next section \ref{2} we express in subsection \ref{2-5} the coefficients via loop counting that enables to perform some explicit computations of the corresponding GF's. In  subsection \ref{3} we recall some properties of circulants and then in the next subsection \ref{proof-HN} we prove a combinatorial identity from which an interesting 
 conjecture in \cite{nagy} follows.

\section{The main theorem}\label{1}

We use (ordinary) generating functions in the ring of formal power series $\QQ[\![t]\!]$  to enumerate DLP's. For convenience {\em hereinafter we consider even-lengthed  DLP's in $\ZZ_+ \times \ZZ^d$
only\footnote{In our setting (as DLP do not hit $T$-axis at the odd moments) the following equalities hold: $\text{\footnotesize{the number of odd-lengthed DLP's of the length}}\, 2k+1 = 2^d\times \text{\footnotesize{the number of even-lengthed DLP's of the length}}\, 2k, k=0,1,\dots$.} and assume that the degree of the formal parameter is 
always equal to the half length of a path. 

Hereinafter, we use abbreviation GF for the {\bf generating function in $\QQ[\![t]\!]$ of the number of objects involved}}.

Recall that a {\em multisection} of a power series or a Laurent series is a new power (or a Laurent) series composed of equally spaced terms extracted unaltered from the original series. Formally, if one is given a power series
$G(t) = \sum_{i=0}^\infty g_i t^i$
then its multisection is a power series of the form
$[G(t)]_{q,r}\stackrel{def}{=}\sum_{i=0}^{\infty} g_{qi+r} t^{qi+r}$, 
where $r, q$ are integers, with $0\leq r < q$. In particular, by definition
$G(t)=\sum_{i=0}^{q-1} [G(t)]_{q,i}$.

As is well known\footnote{See, for instance, \cite{Knuth} or Wiki. This receipt is attributed to Thomas Simpson (1751).} a closed expression for the multisection in terms of the function $G(t)$ may be obtained via discrete Fourier transform (DFT) as follows.

Set 
$\omega_q\stackrel{def}{=}  \exp{(-\frac{2\pi i}{q})}$ and let
${\F}_q=\left[\omega_q^{kj}\right]_{k,j=0,\dots,q-1}$ and, respectively,  
${\F}^{-1}_q=\frac1{q}\left[(\omega_q^{-1})^{kj}\right]_{k,j=0,\dots,q-1}$ be the matrices of DFT and, respectively, of the inverse DFT of order $q$. 

Given a function $G(t)$ and a root of unity $\omega_q$ let us 
denote a vector $\overrightarrow{G(t),\omega_q} \stackrel{def}{=}[G(\omega_q^0\cdot t),\dots, G(\omega_q^{q-1}\cdot t)]^{tr}$
(here $(\cdot)^{tr}$ is a transposition of a vector). It
holds:
\begin{equation}\label{di}
[[G(t)]_{q,0},\dots, [G(t)]_{q,q-1}]^{tr} = {\F}^{-1}_q\,\overrightarrow{G(t),\omega_q}, 
\end{equation}
 
Now as multiplication by $t$ results in the shift of the coefficients of any GF by one position to the right the following identity for the multisections of the GF $t\cdot G(t)$ should hold
\begin{equation}\label{shift}
[t\cdot G(t)]_{q,i} = t[G(t)]_{q,(i-1)\pmod{q}},\; i=0,1,\dots,q-1.
\end{equation}

To formulate the main theorem below let us recall the notions of an {\em excursion} and an {\em escaping DLP}.

An {\em  excursion} is a DLP  whose start  and end points are fixed. Usually the starting point is the origin.
Let call $T$-{\em  excursion} a DLP of even length that starts  and ends    somewhere on 
the $T$-axis.  By definition, a {\em primitive} $T$-excursion should reach $T$-axis at the ends only.

By definition an {\em escaping DLP} should start somewhere at an even point on $T$-axis and never return to it further.

If some points of the lattice are forbidden then by definition the excursions and escaping DLP's should avoid them so that, for instance, primitive $T$-excursions should start and finish at admissible points and an escaping DLP should start from an admissible point. 

\begin{remark}
{\bf A crucial observation.} Some sets of objects involved, say, DLP's, may be not affected by prohibitions and it follows that their numbers, for instance, may be computed on unrestricted lattice as if there were no forbidden points. For example,
just by definition GF's of escaping DLP's in a restricted lattice starting at any admissible point on $T$-axis are equal and may be computed on unrestricted lattice. By the same reasoning the number of primitive $T$-excursions between two admissible points on $T$-axis may either be computed on unrestricted lattice etc. 

Especially we shall use the fact that if the whole multisection of a GF involved is not affected by forbidding  a set of lattice points then {\bf we can compute it as a multisection of the original GF calculated on the unrestricted lattice}. 
\end{remark}

{\em Hereinafter the left superscript in the notation of GF's will indicate the dimension of the underlying ``space'' lattice $\ZZ^d$.}

Let us denote $^dE(t)$ the GF of starting from the origin primitive $T$-excursions on unrestricted lattice $\ZZ_+ \times \ZZ^d$. 
Respectively, let us denote $^dE^{\infty}(t)$ the GF of starting from the origin escaping DLP's on unrestricted lattice $\ZZ_+ \times \ZZ^d$.

It may be helpful here to recall some well known 
facts  
for the $(1+1)$-dimensional DLP's. In this case $T$-excursions correspond either to {\em Dyck}-paths (positive or negative, that are located in the upper or, respectively, the lower halfplane) or concatenation of several Dyck paths. Primitive excursions correspond to Dyck paths (positive or negative) that touch $T$-axis at the ends only.

Recall that the number of Dyck paths 
of length $2k$ is equal to the $k$-th Catalan number 
$B_k = \frac{1}{k+1} \binom{2k}{k}$ 
so  that their GF is given by $CAT(t)=\sum_{k=0}^{\infty} B_k t^k=\frac{1-\sqrt{1-4t}}{2t}$. Note that the 
number of {\em primitive $T$-excursions} of length $2l,\, l> 0$ (positive and negative)  that reach $T$-axis at the ends of path  only equals twice the number of Dyck paths (positive and negative) of length $2(l-1)$ and thus
\begin{equation}\label{excur}
^1E(t) = 2t\cdot CAT(t) = 1-\sqrt{1-4t}. 
\end{equation}

The number of escaping restricted DLP's 
starting from 
any admissible point $(2r,{0})$  on $T$-axis
is equal to the corresponding GF of {\em unrestricted} escaping DLP's starting from the origin, 
and thus for the $(1+1)$-dimensional  case, see, e.g. \cite[Ch. 3, sect.3, Lemma 1]{feller}, it equals:
\begin{equation}\label{escape1}
^1{E}^{\infty}(t)=\frac{1}{\sqrt{1-4t}}.
\end{equation} 
Note that the approach used in \cite{feller} is not directly applicable to counting DLP's in higher dimensions\footnote{As it uses, for instance, the reflection principle.}, but the projection along $T$-axis onto ``space'' part $\ZZ^d$ induces bijections between DLP's and walks in $\ZZ^d$. In particular, excursions correspond to starting from the origin closed walks or loops, primitive $T$-excursions correspond to simple loops (starting from the origin irreducible loops --- see subsection \ref{2-5}) etc. and we may apply loop counting on integer lattice $\ZZ^d$ to compute the numbers of DLP's involved. Some difficulties of course arise as primitive $T$-excursions, for instance, not only should start at an admissible point but should also end at an admissible point either so that some of their lengths should be prohibited but under periodicity limitations certain multisections are not affected by forbidding some periodic sets on $T$-axis. In any case loop counting approach enables to obtain formula for all dimensions in a universal manner.
We'll discuss this technique in \ref{2-5}. 


Denote $^d{P}^r(A,t)$ 
the corresponding 
GF of the number of DLP's of even length starting from an admissible point $(2r,0)$ that may reach $T$-axis at points of an admissible set $A$ only.

Note that if an admissible set $A$ is periodic then by symmetry for any admissible point $(2r \geq 0,0)$ there is a bijection between the sets $^d{\cal P}^{r}(\mathbf{A},t)$ and $^d{\cal P}^{r\pmod{t_A)}}(\mathbf{A},t)$ and hence, the corresponding GF's are equal:
\begin{equation}\label{period}
^d{P}^r(\mathbf{A},t) = ^d{P}^{r\pmod{t_A)}}(\mathbf{A},t).
\end{equation}

Let $j,q\in \{a_0 (=0),\dots,a_k\} \subseteq \{0,1,\dots,t_A-1\}$ be two (possibly equal) elements in a periodic set $\mathbf{A}=(\{a_0,\dots,a_k\}, t_A)$. Let us define {\em shift} operation as follows: 

$$
Sh(j,q)\stackrel{def}{=}\left\{
\begin{array}{cl}
q-j,& \mbox{if } j  \leq q;\\
n + j - q,& \mbox{if } j > q.
\end{array}
\right.
$$

Informally, $Sh(j,q)$ is equal to the distance from the point $(2j,0)$  {\em to the nearest to it to the right on $T$-axis} point $(2y,0),\, y\pmod{t_A}=q, \, y=j+S(j,q)$. 

Another interpretation of $Sh(\cdot,\cdot)$ is as follows. Let $D$ be an oriented graph of a regular $t_A$-gon with the set of vertices $V(D)=\{0,\dots,t_A-1\}$ whose arcs are oriented clockwise. Then  $Sh(j,q),\, j,q\in V(D)$ is equal to graph distance $dist_D(i,j)$, i.e. minimal number of arcs in the path in $D$ between vertices $j$ and $q$. 

Let $A$ be a periodic set and let $^dE^u(\mathbf{A},t)$ the GF of starting from $(2u,0)$ primitive $T$-excursions on restricted lattice $\ZZ_+ \times \ZZ^d$. Assume that the multisection $H=[^dE^u(\mathbf{A},t)]_{t_A,q},\ 0\leq q,u < t_A$   {\bf is not affected by introducing forbidden points} then by definition of the shift operation it holds:
\begin{equation}\label{prim-excur}
[^dE^u(\mathbf{A},t)]_{t_A,q}=[^dE(t)]_{t_A,Sh(u,q)},
\end{equation}
and the multisection $H$ may be computed on the unrestricted lattice.

\begin{theorem}\label{th1} 
Let $\mathbf{A} = (A=\{a_0(=0),\dots,a_k\}, t_A)$  be periodic. 
The set of generating functions $^dP^r(\mathbf{A}
,t), r\in A$ of DLP's on $\ZZ_+\times\ZZ^d$ 
satisfies the following system  of linear equations in the ring of formal power series $\QQ[\![t]\!]$:

\begin{equation}\label{system}
^d{P}^r(\mathbf{A},t) - \sum_{q\in A} {[^d{E}(t)]_{(t_A,\,Sh(r,q)})}\, ^d{P}^q(\mathbf{A},t) = ^d{E}^{\infty}(t),\, r\in A.
\end{equation}

Moreover \eqref{system} 
has a unique solution. 
\end{theorem}

\begin{proof}
 
Let $A$ be a set of admissible points (not necessarily periodic).
Let $^d{E}^{(r_1,r_2)}(A,t)\, \, r_1<r_2$ denote the GF of primitive $T$-excursions starting from an admissible point $R_1=(2r_1,{0})$ and ending at an admissible point $R_2=(2r_2,0)$. Evidently $$^d{E}^{(r_1,r_2)}(A,t) = t^{\text{\footnotesize{number of primitive excursions from A to B}}}.$$ 

Any DLP starting from any admissible point either touches $T$-axis or not. In the former case it can be uniquely decomposed into a primitive $T$-excursion from the starting point to some admissible point $(2r>0,0)$ and start over. Passing to generating functions we by definition get the following infinite set of equations in the ring of formal power series: 
\begin{equation*}
^d{P}^{a_i}(A,t) - \sum_{j\geq i}
^d{E}^{(a_i,a_j)}(A,t) ^d{P}^{a_j}(A,t) - ^d{E}^{\infty}(t) = 0, \, a_i\in A.
\end{equation*}

Let now $A$ be periodic then using identities from \eqref{period}  
we reduce the unknowns to a finite set $^d{P}^{a_i}(\mathbf{A},t), a_i\in \{a_0(=0),\dots,a_k\}$. Let us partition the set of admissible points on $T$-axis into a disjoint union of residue classes modulo $2t_A: \,\sqcup_{i=0}^k \{a = 2a_i\pmod{2t_A},\, a\in\ZZ_+\}$
After collecting the like summands  we obtain as coefficients of unknowns some {\em series multisection} with step $t_A$ (recall the half length normalization of the formal parameter) of the GF's of the unrestricted primitive $T$-excursions $^d{E}^r(A,t),\, 0\leq r\leq t_A$ that are not affected by the forbidden points and hence using the equalities \eqref{prim-excur} we convert the system into the declared form \eqref{system}. 

To prove the existence and uniqueness of the solution of the linear system \eqref{system} let us check that its determinant  is nonzero\footnote{A bit more generally,
one can show
(see, e.g. discussion in https://math.stackexchange.com/questions/2240050/requirements-on-fields-for-determinants-to-bust-dependence.)
that in any commutative ring the kernel of a square matrix is trivial if and only if its determinant is neither zero nor a zero-divisor, and a square matrix with entries in a commutative ring is invertible iff its determinant is invertible in the ground ring.} . 

At first note that the coefficients of the system \eqref{system} may be computed for 
{\em unrestricted} DLP's. We have already indicated this fact for GF's of escaping DLP's. The same argument holds for $[^d{P}^r(\mathbf{A},t)]_{(t_A,\,Sh(r,q)}$. Moreover via projection the computation of the coefficients involved is related to the famous G. Polya's problem \cite{Polya} on return probability to the origin of a lattice walk on $\ZZ^d$. As return probabilities as well as their arbitrary multisections are holomorphic\footnote{In other words each of GF's involved 
is a series with a non zero radius of convergence, and all the identities involved as well as all such identities that would be listed below may be considered
as identities in the algebra of analytic functions on some open unit disc in $\CC$.} the system matrix tends to identity matrix as $t$ tends to $0$. 

\end{proof}

Note that in periodic case the coefficient's matrix of the system \eqref{system} is a principal submatrix (formed by identic sets of rows and columns whose indices belong to $\mathbf{A}$ ) of a $|t_A|\times |t_A|$ complete {\em circulant matrix} 
\begin{equation}\label{circulant}
\scriptstyle{^d\mathbf{C}[(\{0,1,\dots,t_A-1\},t_A),t]\stackrel{def}{=}\| 1 - ^d{E}[(\{0,1,\dots,t_A-1\},t_A),t]_{(t_A,\,Sh(i,j)}\|_{i,j=0,1,\dots,t_A-1},}.
\end{equation}
 
 {\em Below we use abbreviations $^d\mathbf{C}_{|t_A|}$ and $^d\mathbf{C}(\mathbf{A})$ for the underlying complete circulant matrix  and, respectively, for the system's submatrix of \eqref{system}.}

Let us call a power series in $\QQ[\![t]\!]$ {\em $(\alpha,\beta)$-special} if its exponents form an arithmetic progression  $\beta+i \alpha,\, i=0,1,\dots,\,\alpha \in\ZZ_+,\, \beta\in\{0,\dots,\alpha-1\}$. 
By construction any multisection 
$^dP^r(\mathbf{A},t)_{t_A,q},\, r\in A,\, q\in \{0,\dots,t_A-1\}$ and the  coefficients of the matrix 
$^d\mathbf{C}(\mathbf{A})$ are $(t_A,q)$-special. 
Specifically,  by definition 
$ ^d{E}(\mathbf{A},t)_
 {(t_A, Sh(i,j))}
\, i,j\in\{0,1,\dots,t_A-1\}$ is $[t_A,Sh(i,j)]$-special or, equivalently, $[t_A,dist_{D(i,j)}]$-special power series. 

Note that  the product of $(\alpha,\beta_1)$- special and  $(\alpha,\beta_2)$- special power series is $[\alpha,(\beta_1+\beta_2)\pmod{\alpha}]$-special. 

If we are interested in finding an explicit expression for some concrete multisection $M=^dP^{m}(\mathbf{A},t)_{t_A,q},\, m\in A, q\in \{0,\dots,t_A-1\}$ 
then of course we can find $^d{P}^{a_j}(\mathbf{A},t)$ from the solution of the system \eqref{system} and compute the required multisection by Simpson's formula. 

On the other hand, it turns out that we can reduce the linear system \eqref{system} and rewrite it in terms of multisections of variables $^dP^q(\mathbf{A},t),\, q\in A$ only so that the set of multisections involved in the reduced system (one multisection for each variable including $M$) is uniquely determined after fixing $M$ in any equation of \eqref{system}.  In other words, we can obtain a linear system that is a consequence of \eqref{system} whose variables are $M$ and some other multisections.    

Indeed, let   $\{\underline{^d{P}^r}(\mathbf{A},t)\},\, r\in A$ be the solution of \eqref{system} then $\{\underline{^d{P}^r}(\mathbf{A},t)\}_{t_A,l},$ $l=(dist(m,r)+q)\pmod{t_A},\, r\in A\}$ is a solution of a linear system with the same coefficient's matrix $^d\mathbf{C}(\mathbf{A})$ as \eqref{system} and the constant's vector $^d{E}^{\infty}(t)_{t_A,l},$ $l=(dist(m,r)+q)\pmod{t_A},\, r\in A$. The demonstration consists in taking $(t_A,dist(m,r)+q)\pmod{t_A}$- multisection of the $r$-th equality  $\{\underline{^d{P}^r}(\mathbf{A},t)\} - \sum_{q\in A} {[^d{E}(t)]_{(t_A,\,Sh(r,q)})}\, \{\underline{^d{P}^q}(\mathbf{A},t)\} = ^d{E}^{\infty}(t)$ for $r\in A$. 

Let us call the obtained consequence of \eqref{system} {\em $M$-reduction} of  \eqref{system}.

Let us illustrate  the reduction procedure above for the periodic sets $\mathbf{A_1}$ and $\mathbf{A_2}$. Specifically, we are interested in finding more or less explicit expressions for the GF's of multisections
$^dP^0(\mathbf{A_1},t)_{2,0}$ and $^dP^0(\mathbf{A_2},t)_{4,0}$, 
respectively. In the meaningful terms we want to find the GF's of the number of starting from the origin DLP's of even length in the $\ZZ_+\times \ZZ^d$ that do not visit the points on $T$-axis: $(2i,0),\, i=1 \pmod{2}$ or, respectively, do not visit the points on $T$-axis: $(2i, 0),\, i=2,3\pmod{4}$, 

Let us use theorem \ref{th1} and solve the linear system \eqref{system} that for $\mathbf{A_1}$ consists of the only equation
$^dP^0(\mathbf{A_1},t)  - ^d{E}(t)_{2,\,Sh(0,0)} {^dP^0(\mathbf{A_1},t)} = ^d{E}^{\infty}(t)$ and for its $^dP^0(\mathbf{A_1},t)_{2,0}$-reduction we obtain the following equation $$^dP^0(\mathbf{A_1},t)_{2,0} [1 - ^d{E}(t)_{(2,\,Sh(0,0)}]~=~^d{E}^{\infty}(t)_{2,0}$$ so that 
\begin{equation}\label{A1}
^dP^0(\mathbf{A_1},t)_{2,0}=\frac{^d{E}^{\infty}(t)_{2,0}}{1 - ^d{E}(t)_{2,0}}
\end{equation}
and in particular in view of \eqref{excur} and \eqref{escape1} we get $^1P^0(\mathbf{A_1},t)_{2,0}=\frac{1/2\left(\frac1{\sqrt{1-4t}}+\frac1{\sqrt{1+4t}}\right)}{1/2(\sqrt{1-4t}+\sqrt{1+4t})} = \frac1{\sqrt{1-(4t)^2}}$. 

Computation of the second GF is a bit more tedious. Again we obtain the following linear system 
\begin{eqnarray*}
^dP^0(\mathbf{A_2}, t) -^d{E}(t)_{(4,\,Sh(0,0))} {^dP^0(\mathbf{A_2}, t)} -^d{E}(t)_{(4,\,Sh(0,1))} {^dP^1(\mathbf{A_2}, t)}&=& ^d{E}^{\infty}(t)\\
-^d{E}(t)_{(4,\,Sh(1,0))}{ ^dP^0(\mathbf{A_2}, t)} + [1-^1{E}(t)_{(4,\,Sh(1,1))}] {^1P^1(\mathbf{A_2}, t)}&=& ^d{E}^{\infty}(t)\\
\end{eqnarray*}
whose $M=^d~P^0(\mathbf{A_2},t)_{(4,0)}$-reduction is as follows 
(we fill in the values of $Sh(\cdot,\cdot)$):

\begin{eqnarray*}
^dP^0(\mathbf{A_2}, t)_{(4,0)} -^d{E}(t)_{(4,0)} {^dP^0(\mathbf{A_2}, t)_{(4,0)}} - ^d{E}(t)_{(4,1)} {^dP^1(\mathbf{A_2}, t)_{(4,3)}}&=& ^d{E}^{\infty}(t)_{(4,0)}\\
-^d{E}(t)_{(4,3)} {^dP^0(\mathbf{A_2}, t)_{(4,0)}} + [1-^dE(t)_{(4,0)}] {^dP^1(\mathbf{A_2}, t)_{(4,3)}} &=& ^d{E}^{\infty}(t)_{(4,3)}\\
\end{eqnarray*}

In particular, solving the system above for $d=1$ (for instance, by Cramer's rule and simplifications) we can find\footnote{It may be instructive to perform this computation as final reduction of rather cumbersome expressions seems quite unexpected (at least for the author).}: $^1P^0(\mathbf{A_2}, t)_{(4,0)} = \frac1{\sqrt{1- (4t)^4}}$. 

From the first glance the following conjecture seems to be a bit unexpected\footnote{In fact the conjecture was originally formulated for the coefficients of GF's after analyzing much more convincing statistics.} but nevertheless it has been stated in a somewhat different form by P.~Hajnal and G.~V.~Nagy in \cite[Conjecture 9a]{nagy}. Let $\mathbf{A_k}=(\{0,1,\dots,k\}, 2k)$ be a periodic set. Then the following identity holds: 
$$
^1P^0(\mathbf{A_k}, t)_{2k,0} = \frac1{\sqrt{1- (4t)^{2k}}}. 
$$
We call this statement {\em HN-conjecture} and prove it in the next section.

In fact we prove a more general statement from which the {\em HN-conjecture} follows as a consequence.

Specifically below we 
prove that Cramer's solution for $^dP^0(\mathbf{A_k}, t)_{2k,0}$ of the $^dP^0(\mathbf{A_k}, t)_{2k,0}$-reduction of \eqref{system} gives the desired result for $d=1$. Namely, it holds
\begin{equation}\label{HN-conj}
[\mbox{Cramer's determinants ratio for } ^1P^0(\mathbf{A_k}, t)_{2k,0}] = \frac1{det( ^1\mathbf{C}_{2k})} = \frac1{\sqrt{1- (4t)^{2k}}}. 
\end{equation}

\section{Loops. Circulants. Proof of 
HN-conjecture}\label{2}

\subsection{Loops}\label{2-5}

At first we show how to express all the coefficients and the free term of the system \eqref{system} in all dimensions via loop counting on integer lattice $\ZZ^d$.

\emph{A loop is a closed standard lattice walk on a lattice $\ZZ^d$ that {\bf begins
and ends at the origin}}, i.e. a loop is a finite sequence of vectors from $\{\pm 1\}^d$ (steps)  that starts from the origin and sums to zero vector (as a sum of vectors). It may be convenient to regard a loop as a word in the alphabet $\{\pm 1\}^d$.
{\em A loop is {\bf simple} if no non empty proper prefix of it sums to zero as a sequence of vectors}.  Equivalently, {\em a loop is simple if it is not a concatenation (as words) of two nontrivial loops}.
Here we should postulate that the walk of
length zero --- the \emph{trivial} loop--- is not simple.

Obviously, a loop 
has an even length. Let $^d{L}_k$ and, respectively,  $^d{SL}_k$ be the numbers\footnote{Recall that by our convention the left superscript refers to the dimension of the ''space'' lattice involved.} of loops and simple loops of
length $2k$. An easy counting shows: 
\begin{multline*}
^1{L}_0=;\dots,=^d{L}_0=1,\,^1{SL}_0=,\dots,=^d{SL}_0=0;\\ ^1{L}_1=^1{SL}_1=2;\\ ^2{L}_1=^2{SL}_1=4;\, ^2{L}_2=36,\, ^2{SL}_2=20;
\end{multline*}

Let $^dL(t)$ and, respectively, $^d{SL}(t)$ be the corresponding GF's of loops and simple loops.
The identity below (see, e.g. \cite{Novak}) is well known:
$$^dL(t) = ^dL(t)\cdot ^d{SL}(t)+1.$$ 

The following proposition follows directly from the definitions.

\begin{proposition} The projection along $T$-axis establishes bijection between starting from the origin primitive $T$-excursions in $\ZZ\times \ZZ^d$ and starting from the origin simple loops in $\ZZ^d$.
\end{proposition}

Hence, the GF of starting from the origin $T$-excursions on the unrestricted lattice
may be calculated as follows: 
\begin{equation}\label{s-loops}
^dE(t) =^d{SL}(t) = 1 - \frac1{^dL(t)}.
\end{equation} 
For $d=1$ we get $^1L_k = {\binom {2k} k}$ so that $^1L(t) = {\sum_{k=0}^{\infty}} ^1L_k t^k = \frac1{\sqrt{1-4t}}$ and from \eqref{s-loops} it is exactly the GF of primitive excursions already obtained in \eqref{excur}:
$^1E(t) = ^1L(t)  = 1-\sqrt{1-4t}$.

{\em Note that in general the GF obtained is not equal to the corresponding GF of primitive excursion on a restricted lattice but all its multisections involved in the coefficients matrix of the \eqref{system} are not affected by prohibitions and may be calculated as corresponding multisections of the GF defined in \eqref{s-loops}.
}

Computing the free term of \eqref{system} is the same as computing the GF of escaping paths 
$^d{E}^{\infty}(t) = \sum_{k=0}^{\infty} [^d{E}^{\infty}]_k t^k$ and we now show how to reduce it (for all $d$) to loop counting. Namely, any {\em non escaping DLP} of the length $2k$ may be uniquely decomposed as a concatenation of a starting from the origin 
primitive excursion of the length $2i,\,i=1,\dots,k$ followed by starting from $(2k,0)$ an arbitrary DLP's of the length $2(k-i)$ so that 
$[^d{E}^{\infty}]_k = (2d)^{2k} - \sum_{i=1}^k {^d{SL}}_i (2d)^{2(k-i)}$ and for the GF we get:
\begin{equation}\label{free-term}
^d{E}^{\infty}(t) = \frac1{1-4d^2t} - ^d{SL}(t)\cdot \frac1{1-4d^2t} = \frac1{^dL(t) (1-4d^2t)}.
\end{equation}

{\em And as above all the multisections of GF of starting from an admissible point escaping DLP's may be calculated on unrestricted lattice.}

For $d=1$ we get already obtained in \eqref{escape1} the GF  of escaping DLP's:    $^1{E}^{\infty}(t) = \frac{\sqrt{1-4t}}{1-4t} = \frac1{\sqrt{1-4t}}$.

Thus using \eqref{s-loops}, \eqref{free-term}, 
and \eqref{prim-excur} we may directly rewrite \eqref{system} for an arbitrary periodic set in all dimensions in an ``invariant'' form via GF's of loops and/or simple loops. As an example \eqref{A1} may be rewritten as follows
\begin{equation}\label{A2}
\scriptstyle{
^dP^0(\mathbf{A_1},t)_{2,0}=\frac{
\frac1{^dL(t) (1-4d^2t)}+\frac1{^dL(-t) (1+4d^2t)}}{\frac1{^dL(t)}+\frac1{^dL(-t)}} =
\frac1{1-16d^4t^2}\left[1 - 4d^2t \frac{^dL(t)-^dL(-t)}{^dL(t)-^dL(-t)}\right].}
\end{equation}

Moreover as there are some explicit enough analytic expressions for loops and/or simple loops we can continue simplifications of \eqref{system}. For instance,
it has been already indicated in the original Polya's paper \cite[S. 160]{Polya} that $^2L(t)$ may be expressed
via the elliptic integral of the first kind. Namely, the \emph{complete elliptic integral of the first kind} with \emph{elliptic modulus} $q$ in the open unit disc is defined by
\begin{equation}
K(q) = \int_{0}^1 \frac{\mathrm{d} y}{\sqrt{(1-y^2)(1-q^2y^2)}},
\end{equation}. 
and has
series expansions around the origin given by \cite[17.3.11-12]{abramovitz}
\begin{equation}\label{eq:Kexpansion}
K(q) = \frac{\pi}{2} \sum_{k=0}^\infty {\binom{2k}{k}}^2 \left(\frac{q}{4}\right)^{2k}, \qquad 
\end{equation}
that gives an analytic expression for $^2L(t)=2/\pi K(4 \sqrt{t})$ and the expression \eqref{A2} for $d=2$ may be rewritten as follows
\begin{multline}\label{A3}
\scriptstyle{^2P^0(\mathbf{A_1},t)_{2,0} =
\frac1{1-256t^2}\left[1 - 16t \frac{K(4\sqrt{t})-K(4\sqrt{-t})}{K(4\sqrt{t})+K(4\sqrt{-t})}\right] =}\\
\scriptstyle{= 1+192 t^{2}+45056 t^{4}+10979328 t^{6}+2716942336 t^{8}+677907697664 t^{10}+170013263888384 t^{12}+\dots
}
\end{multline}

Analogously somewhat more complex analytic expressions for $^d{SL}(t)$ are also known and may be used. 
Namely, 
it is shown 
in \cite{Novak} that
\begin{equation}\label{SL} 
		^d{SL}(t) 
		= \int_0^{\infty} I_0\bigg{(}\frac{tu}{d}\bigg{)}^d e^{-t} du,
	\end{equation}
where $I_0(z)$ is the {\em modified Bessel function} (see, e.g.  \cite[Chapter 4]{AAR}) that is a solution of the \emph{modified Bessel differential equation} and
admits
both a series representation,

	\begin{equation*}
		I_\alpha(z) = \sum_{k=0}^{\infty} \frac{(\frac{z}{2})^{2k+\alpha}}{k!\Gamma(k+\alpha+1)},
	\end{equation*}
	
\noindent
and an integral representation,

	\begin{equation*}
		I_\alpha(z) = \frac{(\frac{z}{2})^\alpha}{\sqrt{\pi}\Gamma(\alpha+\frac{1}{2})}\int_0^\pi 
		e^{(\cos \theta)z}(\sin \theta)^{2\alpha} d\theta.
	\end{equation*}
	

It may be relevant to note here that in combinatorial setting the complexity of loop counting is usually expressed in terms of the analytic class in which  the GF,s involved falls (e.g., rational, algebraic, transcendental, $D$-finite etc). And from this viewpoint the complexity grows drastically when dimension grows. Namely, $^1L(t)$ is not rational but algebraic; $^2L(t)$ is not algebraic but transcendental and $D$-finite; $^3L(t)$ according to \cite{dfinite} is not  $D$-finite 

\subsubsection{Circulants}\label{3}

Our proof  of HN-conjecture relies on the properties of the circulant matrices as the coefficient matrix $^d\mathbf{C}(\mathbf{A})$ of the system \eqref{system} is a principal submatrix of the circulant matrix $^d\mathbf{C}_{|t_A|}$.

Let us recall some standard facts about {\em circulants} that may be found e.g., in  \cite{circulant} or even in Wikipedia. 

A matrix of the form 
$$
C\stackrel{def}{=}
\begin{bmatrix}
c_0&c_{1}&\dots&c_{n-1}\\
c_{n-1}&c_0&\dots&c_{n-2}\\
\vdots&\vdots&\ddots&\vdots\\
c_{1}&c_{2}&\dots&c_0\\
\end{bmatrix}
$$
is called a {\em circulant matrix}. Let us denote $c^{row}=(c_0,\dots,c_{n-1})$  its first row vector and, respectively, let us denote $c^{col}=(c_0,c_{n-1},\dots,c_{1})^{tr}$ its first column vector (recall that $(\cdot)^{tr}$ is a transpose of a vector).

A circulant matrix $C$  may be fully specified by the first row $c^{row}$ (a choice in MAPLE) or, respectively, by the first column  $c^{col}$  (a choice in Wikipedia).  The remaining lines (rows or columns) of $C$ are cyclic permutations of the corresponding lines with offset equal to the column or row index minus $1$.

$c^{row}$ and $c^{col}$ are related by the following identity:
\begin{equation}\label{row2col}
c^{col}=W (c^{row})^{tr},
\end{equation}
where  
$$
W\stackrel{def}{=}
\begin{bmatrix}
\begin{array}{ccccc}
1&0&\dots&0&0\\
0&0&\dots&0&1\\
0&0&\dots&1&0\\
\vdots&\vdots&\ddots&\vdots\vdots\\
0&1&\dots&0&0\\
\end{array}
\end{bmatrix}
$$
Note that the following equality holds
\begin{equation}\label{W}
W=\frac1{n} {{\F}_n}^2,
\end{equation}
here ${\F}_n$ is a $n\times n$ matrix of DFT.
In particular, $W$ is a permutation matrix that commutes with ${\F}_n$.

{\em Hereinafter we assume that a circulant is specified by the first row and denote $CM(c^{row})$ a circulant matrix whose first row equals $c^{row}$}.

Let us denote
the first row and the first column of the underlying complete circulant matrix $^d\mathbf{C}_{|t_A|}$ of the system \eqref{system} 
 by  $\widetilde{c^{row}}$ and $\widetilde{c^{col}}$, respectively, so that it holds $^d\mathbf{C}_{|t_A|}=CM(\widetilde{c^{row}})$.

By definition we have
\begin{multline}\label{firstrow}
\scriptstyle{
\widetilde{c^{row}}  = (c_0,c_1,\dots,c_{|t_A|-1}) 
 \stackrel{def}{=}[1-^dE_{|t_A|,Sh(0,0)},-^dE_{|t_A|,Sh(0,1)},\dots,-^dE_{|t_A|,Sh(0,|t_A|-1)}]=}\\
\scriptstyle{ [1-(1-\frac1{^dL(t)_{|t_A|,Sh(0,0)}},-(1-\frac1{^dL(t)_{|t_A|,Sh(0,1)}}),\dots, -(1-\frac1{^dL(t)_{|t_A|,Sh(0,|t_A|-1)}}]}
\end{multline} 
and
\begin{equation}\label{firstcol}
\widetilde{c^{col}}  = (c_0,c_{|t_A|-1},\dots,c_2)^{tr}
\end{equation}



\noindent 
Let us denote\footnote{The subscript indicates the dimension.} $CRC_{|t_A|}\stackrel{def}{=}\overrightarrow{\frac1{^dL(t)},\omega_{|t_A|}}= \scriptstyle{ \big[\frac1{^dL(t\cdot 1)},\frac1{^dL(t\cdot \exp{(-\frac{2\pi i 1}{|t_A|})})},\dots,\frac1{^dL(t\cdot \exp{(-\frac{2\pi i (|t_A|-1)}{|t_A|})})}\big]^{tr}}$
and
let $\Lambda=(\lambda_1,\dots,\lambda_{|t_A|})^{tr}$ be the vector of eigenvalues 
of the circulant matrix $^d\mathbf{C}_{|t_A|}$. 
It is well known that $\Lambda$ is equal to
the product of ${\F}_{|t_A|}$ by the first column of $^d\mathbf{C}_{|t_A|}$:
\begin{multline}\label{eigen}
\scriptstyle{\Lambda = {\F}_{|t_A|} \widetilde{c^{col}} = {\F}_{|t_A|} W (\widetilde {c^{row}})^{tr} = W {\F}_{|t_A|} (\widetilde{c^{row}})^{tr} =
W {\F}_{|t_A|} (\underbrace{(1,0,\dots,0)}_{|t_A|})^{tr}-W(\underbrace{(1,1,\dots,1)}_{|t_A|})^{tr}+}\\ 
\scriptstyle{
W {\F}_{|t_A|} ({\F}_{|t_A|})^{-1}\big[\frac1{^dL(t\cdot 1)},\frac1{^dL(t\cdot \exp{(-\frac{2\pi i 1}{|t_A|})})},\dots,\frac1{^dL(t\cdot \exp{(-\frac{2\pi i (|t_A|-1)}{|t_A|})})}\big]^{tr}= W\cdot CRC_{|t_A|}},
\end{multline}
Hence, in the Fourier eigenbasis we obtain 
 the following representation: 
 \begin{equation}\label{fourier-repr}
 ^d\mathbf{C}_{|t_A|} = ({\F}_{|t_A|})^{-1} diag[W\cdot CRC_{|t_A|}] {\F}_{|t_A|}, 
\end{equation} 
(here $diag(\cdot)$ denotes a diagonal matrix built from a vector) and it follows directly from \eqref{fourier-repr} that the first row of $^d\mathbf{C}_{|t_A|}$ equals $({\F}_{|t_A|})^{-1} CRC_{|t_A|}$ so that $^d\mathbf{C}_{|t_A|}=CM[({\F}_{|t_A|})^{-1} CRC_{|t_A|}]$. Analogously we conclude that the following equality holds for the inverse circulant matrix $(^d\mathbf{C}_{|t_A|})^{-1} = CM[{\F}_{|t_A|}^{-1} \frac1{CRC_{|t_A|}}]$, 
where by definition for a vector $\stackrel{\to}{v}=(v_1,\dots,v_n)$ with non zero components we set   $\frac1{\stackrel{\to}{v}}\stackrel{def}{=}(\frac1{v_1},\dots,\frac1{v_n})$.

Thus given a periodic set of admissible points $\mathbf{A} = (A=\{a_0(=0),\dots,a_k\}, t_A)$  on $T$-axis we can 
substitute the elements of the $|A|\times|A|$ circulant matrix $^d\mathbf{C}(\mathbf{A})$ and, respectively, free term with their expressions through loops and/or simple loops \eqref{s-loops} and \eqref{free-term} 
and obtain an equivalent form of the system \eqref{system}: 

\begin{equation}\label{system-loops}
\big[\widehat{^d\mathbf{C} (\mathbf{A})}\big] \stackrel{\to}{P}  = 
\overrightarrow{^d{E}^{\infty}(t)},
\end{equation}
here $\big[\widehat{^d\mathbf{C}(\mathbf{A})}\big]$ is a $|A|\times|A|$ principal submatrix of complete circulant matrix $CM[({\F}_{|t_A|})^{-1} CRC_{|t_A|}]$ 
written in the Fourier eigenbasis,
whose rows and columns are numbered by the set $A$, $\stackrel{\to}{P}$ is a vector $(^d{P}^r(\mathbf{A},t))^{tr}, r\in A$, and $\overrightarrow{^d{E}^{\infty}(t)}$ is an 
$|A|$-dimensional vector all of whose components are equal to $^d{E}^{\infty}(t)$.

In particular for the periodic set $\mathbf{A_k}$ the $k\times k$ matrix $\widehat{^d\mathbf{C} (\mathbf{A_k})}$ of the system \eqref{system-loops} is equal to the left upper quarter of the complete circulant matrix $^d\mathbf{C}_{2k}$, whose first row is equal to $\widehat{c^{row}}\stackrel{def}{=} (\widehat{c_0},\widehat{c_1},\dots,\widehat{c_{2k-1}}) = ({\F}_{2k})^{-1}\overrightarrow{\frac1{^dL(t)},\omega_{2k}} = ({\F}_{2k})^{-1} CRC_{2k}$. 

Additionally let us define a vector $b^{row}=(b_0,b_1,\dots,b_{2k-1})=({\F}_{2k})^{-1}\overrightarrow{[\frac1{^dL(t)(1-4d^2t)},\omega_{2k}]}$  and the corresponding circulant 
$^d\mathbf{B}_{2k}\stackrel{def}{=}CM(b^{row})$.
Performing $^d{P}^0(\mathbf{A_k},t)_{(2k,0)}$-reduction of the system involved related to $\mathbf{A_k}$ we get the following linear system
\begin{equation}\label{system1}
\begin{bmatrix}
\widehat{c_0}&\widehat{c_1}&\dots&\widehat{c_k}\\
\widehat{c_{2k-1}}&\widehat{c_0}&\dots&\widehat{c_{k-1}}\\
\vdots&\vdots&\ddots&\vdots\\
\widehat{c_{k+1}}&\widehat{c_{k+2}}&\dots&\widehat{c_{0}}\\
\end{bmatrix}
\begin{bmatrix}
^d{P}^0(\mathbf{A_k},t)_{(2k,0)}\\
^d{P}^1(\mathbf{A_k},t)_{(2k,2k-1)}\\
\vdots\\
^d{P}^{k-1}(\mathbf{A_k},t)_{(2k,k+1)}\\
\end{bmatrix}
=
\begin{bmatrix}
^d{E}^{\infty}(t)_{(2k,0)}\\
^d{E}^{\infty}(t)_{(2k,2k-1)}\\
\vdots\\
^d{E}^{\infty}(t)_{(2k,k+1)}\\
\end{bmatrix}
=\begin{bmatrix}
b_0\\
\\b_{2k-1}\\
\vdots\\
b_{k+1}\\
\end{bmatrix}\\
\end{equation}

Evidently any even circulant matrix of the order $2k$ may be partitioned into four $k\times k$ ``quarter'' matrices so that $^d\mathbf{C}_{2k}=
\big[
\begin{smallmatrix}
{C_1}&{C_2}\\
{C_2}&{C_1}
\end{smallmatrix}
\big]$
and $^d\mathbf{B}_{2k}=
\big[
\begin{smallmatrix}
{B_1}&{B_2}\\
{B_2}&{B_1}
\end{smallmatrix}
\big],$
here $C_i, B_i,\, i=1,2$ are corresponding $k\times k$ submatrices. By Cramer's rule we   get $^dP^0(\mathbf{A_k}, t)_{2k,0}=\frac{det(\widehat{C_1})}{det(C_1)}$, where 
$\widehat{C_1}$ differs from $C_1$ by the first column that is equal to a constant vector of the system \eqref{system1}, i.e. by definition to the first column of $B_1$. 

Now we prove that $det(\widehat{C_1})=det(B_1)$.

At first let us prove that
%
$b_i - \widehat{c_i} = 4d^2tb_{(i-1)~\pmod{2k}},\, i=0,1,\dots,2k-1$, i.e. the difference of the $i$-th $(i=1,\dots,2k-1)$ components of the vectors $b^{row}$ and $\widehat{c^{row}}$ is equal to a multiple 
of the $(i-1)~\pmod{2k}$-th component of $b^{row}$. 

Indeed,
\begin{multline}\label{b-c}
b_i - \widehat{c_i} = 
\Biggl[
({\F}_{2k})^{-1}
\left(
\overrightarrow{\frac1{^dL(t)(1-4d^2t)},\omega_{2k}} - \overrightarrow{\frac1{^dL(t)},\omega_{2k}}
\right)
\Biggr]_i\\
\Biggl[ 
4d^2 ({\F}_{2k})^{-1}\overrightarrow{\frac{t}{^dL(t)(1-4d^2t)},\omega_{2k}}
{\Biggr]}_i 
= 4d^2t b_{(i-1)~\pmod{2k}},\, i=0,1,\dots,2k-1,
\end{multline}
the last equality follows from \eqref{shift}.

\begin{lemma} 
Let $\tilde C$ be any square submatrix of $^d\mathbf{C}_{2k}$ formed by consecutive columns and let $\tilde B$ be the corresponding submatrix in $^d\mathbf{B}_{2k}$. If we substitute the first column of $\tilde B$ for the first column of $\tilde C$ then\footnote{In fact lemma is valid  for any pair of circulant matrices $Q$ and 
$Q\cdot ({\F})_q^{-1} diag\left[\overrightarrow{\frac1{1-(const)\dot t},\omega_q}\right] \F_q$ ($q$ is the order of $Q$).} 
the determinant of the newly formed matrix equals $det(\tilde B)$. 

\end{lemma}

\begin{proof}
It is enough to consider the case of the submatrix $\tilde C=\widehat{C_1}$ as the general statement may be proven analogously. Indeed, it follows from \eqref{b-c} that all matrices in the following sequence, in which each neighboring pair differs by one column operation\footnote{The first matrix in the sequence is $\widehat{C_1}$ i.e. $C_1$ with  altered first column and the last one  is ${B_1}$ and each next matrix in the sequence differs from the previous by a column operation not changing its determinant.} evidently have equal determinants:
\begin{multline*}
\scriptstyle{
\begin{bmatrix}
b_0&\widehat{c_1}&\dots&\widehat{c_{k-1}}\\
b_{2k-1}&\widehat{c_0}&\dots&\widehat{c_{k-2}}\\
\vdots&\vdots&\ddots&\vdots\\
b_{k+1}&\widehat{c_{k+2}}&\dots&\widehat{c_0}\\
\end{bmatrix}
}
\to
\begin{bmatrix}
b_0&b_1-4d^2t b_0&\dots&\widehat{c_{k-1}}\\
b_{2k-1}&b_0-4d^2t b_{2k-1}&\dots&\widehat{c_{k-2}}\\
\vdots&\vdots&\ddots&\vdots\\
b_{k+1}&b_{k+2}-4d^2t b_{k+1}&\dots&\widehat{c_0}\\
\end{bmatrix}
\to \dots \to\\
\begin{bmatrix}
b_0&b_1&\dots&b_{k-2}&b_{k-1}-4d^2t b_{k-2}\\
b_{2k-1}&b_0&\dots&b_{k-3}&b_{k-2}-4d^2t b_{k-3}\\
\vdots&\vdots&\ddots&\vdots&\vdots\\
b_{k}&b_{k+1}&\dots&b_{0}&b_1-4d^2t b_{0}\\
b_{k+1}&b_{k+2}&\dots&b_{2k-1}&b_0-4d^2t b_{2k-1}\\
\end{bmatrix}
\to \mathbf{B_1}.
\end{multline*}

\end{proof}

Thus the following equality holds
\begin{equation}\label{generHN}
^dP^0(\mathbf{A_k}, t)_{2k,0}=\frac{det({B_1})}{det(C_1)}.
\end{equation}

\subsubsection{Proof HN-conjecture}\label{proof-HN}

Note that incidentally for $d=1$ (as $\frac1{^1L(t) (1-4t)}=\frac{\sqrt{1-4t}}{1-4t}=^1L(t)$) it holds 
$$
^d\mathbf{B}_{2k} = (^d\mathbf{C}_{2k})^{-1}
$$

Finally, let us define a block $2k\times 2k$ matrix  $\mathbf{D}_{2k} =
 \big[
\begin{smallmatrix}
{B_1}&{B_2}\\
0&({C_1})^{-1}
\end{smallmatrix}
\big]$ (here all zero's block is denoted by $0$).
It holds 
$\mathbf{D}_{2k} ^1\mathbf{C}_{2k} =
 \big[
\begin{smallmatrix}
{B_1}&{B_2}\\
0&({C_1})^{-1}
\end{smallmatrix}
 \big]\,
 \big[
\begin{smallmatrix}
{C_1}&{C_2}\\
{C_2}&{C_1}
\end{smallmatrix}
 \big]
=
 \big[
\begin{smallmatrix}
{Id_k}& 0 \\
{(C_1)^{-1}} C_2&{Id_k}
\end{smallmatrix}
 \big]
 $ (here ${Id_k}$ stands for identity $k\times k$ matrix). Equating determinants on both sides we derive the following identity 
 \begin{equation}\label{deter}
\frac{det (B_1)} {det (C_1)} = det [(^1\mathbf{C}_{2k})^{-1}].
\end{equation}
Note that \eqref{deter} holds for any $2k\times 2k$ circulant matrix with a non degenerate left upper quarter matrix.

Thus HN-conjecture follows as
\begin{equation}\label{det-1}
det (^1\mathbf{C}_{2k})^{-1} = \scriptstyle{
\frac1{\sqrt{1-4t}} \frac1{\sqrt{1-4t\exp{(-\frac{2\pi i1}{2k})}}}, \dots, \frac1{\sqrt{1-4t\exp{(-\frac{2\pi i (2k-1)}{2k})}}} = \frac1{\sqrt{1- (4t)^{2k}}}.
}
\end{equation}

\end{document}